\documentclass[12pt]{article}
\usepackage{graphicx}
\usepackage{amssymb}
\usepackage{amsmath}
\usepackage{cite}
\newcommand{\norm}[1]{\left\Vert#1\right\Vert}

\numberwithin{equation}{section}

\newtheorem{Pa}{Paper}[section]
\newtheorem{Tm}[Pa]{{\bf Theorem}}
\newtheorem{La}[Pa]{{\bf Lemma}}
\newtheorem{Cy}[Pa]{{\bf Corollary}}
\newtheorem{Rk}[Pa]{{\bf Remark}}

\title{The sine kernel, two corresponding operator identities, and random matrices}

\author{ Lev Sakhnovich}

\date{}

\begin{document}

\maketitle

\begin{abstract}
In the present paper, we consider the integral operator, which acts in Hilbert space
and has sine kernel. This operator generates two operator identities and two corresponding canonical differential systems.
We find the asymptotics of the corresponding resolvent and  Hamiltonians. We use both the  method of operator identities  and  the theory of random matrices.
\end{abstract}

\section{Introduction}
Let us consider the operator
\begin{equation}S_{\zeta}f=f(x)-\int_{0}^{\zeta}k(x-y)f(y)dy,\quad \zeta>0,\label{1.1}\end{equation}
where
\begin{equation} k(x)=\frac{\sin (x\pi)}{x\pi}.  \label{1.2}\end{equation}
The operator $S_{\zeta}$ is invertible (see \cite[p. 167]{DIZ}). Hence, we have
\begin{equation}S_{\zeta}^{-1}f=f(x)+\int_{0}^{\zeta}R_{\zeta}(x,y)f(y)dy,\quad f(y){\in}L^{2}(0,\zeta),\label{1.3}\end{equation} where the function $R_{\zeta}(x,y)$ is continuous with respect to the variables $x,y,\zeta.$ The operator $S_{\zeta}$ plays an important role in a number of theoretical and applied problems (e.g., in random matrix theory 
\cite{DIZ, Sak5, TrWi} and in optical problems \cite{Lev}). The operator $S_{\zeta}$ satisfies simultaneously two operator identities and generates two canonical differential systems \cite{Sak1, Sak2, Sak3, Sak4, Sak5}. In the present paper, we investigate the asymptotics
of the kernel $R_{\zeta}(x,y)$ and of the Hamiltonians of the corresponding canonical systems when $\zeta{\to}\infty$ (see Theorem 3.3).
In this case, we use both  the method of operator identities  and the results from theory of random matrices (see \cite{TrWi} and  section 4  of the present paper).
Let us introduce the operators:
\begin{equation}K_{\pm}(\zeta)f=\int_{-1}^{1}k_{\pm}(x,t,\zeta)f(t)dt, \quad f(x){\in}L^{2}(-1,1),
\label{1.4}\end{equation}and
\begin{equation}K(\zeta)f=\int_{-1}^{1}k(x,t,\zeta)f(t)dt, \quad f(x){\in}L^{2}(-1,1),
\label{1.5}\end{equation}
where \begin{equation}k(x,t,\zeta)=\frac{\sin \zeta\pi(x-t)}{\pi(x-t)},\label{1.6}\end{equation} \begin{equation} k_{\pm}(x,t,\zeta)=[\frac{\sin \zeta\pi(x-t)}{\pi(x-t)}{\pm}\frac{\sin \zeta\pi(x+t)}{\pi(x+t)}]/2
\label{1.7}\end{equation}
In  section 4, we  investigate the operators
\begin{equation}S(\zeta,\lambda)=I-{\lambda}K(\zeta)\quad and \quad S_{\pm}(\zeta,\lambda)=I-{\lambda}K_{\pm}(\zeta),\quad 0<\lambda{\leq}1.
\label{1.8}\end{equation}
We note that Fredholm determinants
 \begin{equation}P_{\pm}(\zeta,\lambda)=\det(I-{\lambda}K_{\pm}(\zeta)),\quad P(\zeta,\lambda)=\det(I-{\lambda}K(\zeta))
\label{1.9}\end{equation}
play an important role in the random matrix theory \cite{TrWi}. We found special integral representations
for $P_{\pm}(\zeta,\lambda)$ and for $P(\zeta,\lambda)$. These results may be applied to a number of problems (see Remark 4.8).
\section{Two operator identities and two canonical differential systems}
The operator $S_{\zeta}$, which is defined by formulas \eqref{1.1} and \eqref{1.2}, satisfies the following operator identity (see \cite{Sak1}-\cite{Sak5}):
\begin{equation}(QS_{\zeta}-S_{\zeta}Q)f=-\frac{1}{2i\pi}\int_{0}^{\zeta}[e^{i(x-y)\pi}-e^{-i(x-y)\pi}]f(y)dy,
\label{2.1}\end{equation}
where
\begin{equation}Qf=xf(x).\label{2.2}\end{equation}
The second operator identity has the form \cite{Sak4}:
\begin{equation}(AS_{\zeta}-S_{\zeta}A)f=i\int_{0}^{\zeta}[M(x)+M(y)]f(y)dy;\,\label{2.3}\end{equation}
where
\begin{equation}Af=i\int_{0}^{x}f(x)dx,\quad M(x)=\frac{1}{2}-\int_{0}^{x}k(x)dx.\label{2.4}\end{equation}
The operator $S_{\zeta}$  and operator identities \eqref{2.1} and \eqref{2.3} generate two canonical differential systems.The first system is connected with identity \eqref{2.1} and has the form \cite{Sak1}:
\begin{equation}\frac{d}{dx}W_{1}(x,z)=-i\frac{J_1H_{1}(x)}{x-z}W_{1}(x,z),\quad W_{1}(0,z)=I_{2},\label{2.5}\end{equation}
where $J_{1}$ and $H_{1}(x)$ are defined by the relations
\begin{equation}H_{1}(x)=\frac{1}{2\pi}\left(
                           \begin{array}{cc}
                             |q(x)|^2 & q^{2}(x) \\
                             \overline{q^{2}(x)} & |q(x)|^2 \\
                           \end{array}
                         \right);\quad J_{1}=\left(
                                              \begin{array}{cc}
                                                1& 0 \\
                                                0 & -1 \\
                                              \end{array}
                                            \right),
                         \label{2.6}\end{equation}
\begin{equation}q(x)=S_{-}^{-1}e^{ix\pi}.\label{2.7}\end{equation}
The second system is connected with identity \eqref{2.3} and has the form \cite{Sak1}
\begin{equation}\frac{d}{dx}W_{2}(x,z)=izJ_2H_{2}(x)W_{2}(x,z),\quad W_{2}(0,z)=I_{2},\label{2.8}\end{equation}
where $J_{2}$ and $H_{2}(x)$ are defined by the relations
\begin{equation}H_{2}(x)=\frac{1}{2\pi}\left(
                           \begin{array}{cc}
                             q_{1}^{2}(x) & 1/2 \\
                          1/2 &q_{2}^{2}(x)  \\
                           \end{array}
                         \right),\quad J_{2}=\left(
                                              \begin{array}{cc}
                                                0& 1 \\
                                                1 & 0 \\
                                              \end{array}
                                            \right),
                         \label{2.9}\end{equation}
\begin{equation}q_{1}(x)=S_{-}^{-1}1,\quad q_{1}(x)q_{2}(x)=1/2.\label{2.10}\end{equation}
The operator $S_{-}^{-1}  $ is defined by the Krein formula (see \cite{GoKr}Ch.IV and \eqref{1.3}):
\begin{equation}S_{-}^{-1}f=f(x)+\int_{0}^{x}R_{x}(x,y)f(y)dy.\label{2.11}\end{equation}
We note that the operator $S_{\zeta}$ admits the triangular factorization (see \cite{GoKr}Ch.IV):
\begin{equation}S_{\zeta}=S_{-}S_{-}^{\star}.\label{2.12}\end{equation}
\section{Asymptotic behaviour of Hamiltonians $H_{1}(\zeta)$, $H_{2}(\zeta)$ and resolvent kernel $R_{\zeta}(x,y)$
 as $\zeta{\to}\infty$}
 Along  with the operator $S_{\zeta}$ we consider the operator
 \begin{equation}C_{\zeta}f=f(x)-\int_{-\zeta}^{\zeta}k(x-y)f(y)dy,\quad \zeta>0,\label{3.1}\end{equation}
 where the kernel $k(x)$ is defined by \eqref{1.2}.The operator $U_{\zeta}f(x)=f(x-\zeta)$ maps unitary the space
 $L^{2}(-\zeta,\zeta)$ onto $L^{2}(0,2\zeta)$. It is easily to see that
\begin{equation}C_{\zeta}=U_{\zeta}^{-1}S_{2\zeta}U_{\zeta}.\label{3.2}\end{equation}By \eqref{1.3} and
\eqref{3.2} we have
\begin{equation}C_{\zeta}^{-1}f=f(x)+\int_{-\zeta}^{\zeta}Q_{\zeta}(x,y)f(y)dy,\quad f(y){\in}L^{2}(-\zeta,\zeta))\label{3.3}\end{equation}
where
\begin{equation}R_{2\zeta}(x,y)=Q_{\zeta}(x-\zeta,y-\zeta).\label{3.4}\end{equation}
Relation \eqref{3.4} implies that
\begin{equation}R_{2\zeta}(2\zeta,2\zeta)=Q_{\zeta}(\zeta,\zeta),\quad R_{2\zeta}(2\zeta,0)=Q_{\zeta}(\zeta,-\zeta).\label{3.5}\end{equation}
Further we need the following two lemmas.
\begin{La}\label{Lemma 3.1}Let  relations \eqref{3.1} and \eqref{3.3} be fulfilled . Then we have
\begin{equation} Q_{\zeta}(-x,-y)=Q_{\zeta}(x,y).\label{3.6}\end{equation}\end{La}
\emph{Proof.} From \eqref{3.1} and \eqref{3.3} we obtain
\begin{equation}-k(x,y)+Q_{\zeta}(x,y)-\int_{-\zeta}^{\zeta}Q_{\zeta}(x,s)k(s,y)ds=0.\label{3.7}\end{equation}
Relation \eqref{3.7} implies that
\begin{equation}-k(-x,-y)+Q_{\zeta}(-x,-y)-\int_{-\zeta}^{\zeta}Q_{\zeta}(-x,-s)k(-s,-y)ds=0.\label{3.8}\end{equation}

Using equality
\begin{equation} k(-x,-y)=k(x,y).\label{3.9}\end{equation}
and relation \eqref{3.8} we obtain \eqref{3.6}.Lemma is proved.
\begin{La}\label{Lemma 3.2}Let  relations \eqref{3.1} and \eqref{3.3} be fulfilled . Then we have
\begin{equation}q(2\zeta)=e^{i\zeta\pi}r(\zeta),\label{3.10}\end{equation}
where $r(\zeta)$ is defined by the equality
\begin{equation}r(\zeta)=e^{i\zeta\pi}+\int_{-\zeta}^{\zeta}Q_{\zeta}(-\zeta,s)e^{-is\pi}ds.
\label{3.11}\end{equation}
\end{La}
\emph{Proof.} It follows from \eqref{3.11} that
\begin{equation}e^{i\zeta\pi}r(\zeta)=e^{2i\zeta\pi}+\int_{-\zeta}^{\zeta}Q_{\zeta}(-\zeta,-s)e^{i(\zeta+s)\pi}ds.
\label{3.12}\end{equation} Formula \eqref{3.12} can be written in the form:
\begin{equation}e^{i\zeta\pi}r(\zeta)=e^{2i\zeta\pi}+\int_{0}^{2\zeta}Q_{\zeta}(-\zeta,\zeta-t)e^{it\pi}dt.
\label{3.13}\end{equation}Taking into account \eqref{3.4} and \eqref{3.6}, we have
\begin{equation}e^{i\zeta\pi}r(\zeta)=e^{2i\zeta\pi}+\int_{0}^{2\zeta}R_{2\zeta}(2\zeta,t)e^{it\pi}dt.\label{3.14}
\end{equation}The assertion of the lemma follows directly from \eqref{2.7}, \eqref{2.11} and \eqref{3.14}.\\
Let us formulate the main result of the present section:
 \begin{Tm}\label{Theorem 3.1}Let relation \eqref{1.1}, \eqref{1.2} be fulfilled. Then the following asymptotic equalities are valid:\\

1)\begin{equation}R_{2\zeta}(2\zeta,2\zeta){\sim}\pi(\frac{1}{4}u+
\frac{1}{4u}-\sum_{n=1}^{\infty}\frac{c_{2n}}{u^{2n+1}}),\quad u{\to}\infty,\quad (u=2\pi\zeta),\label{3.15}\end{equation}
where $c_2=-\frac{1}{4},\,c_4=-\frac{5}{2}, ...$\\
2)\begin{equation}R_{\zeta}(\zeta,0){\sim}\pi\sum_{n=0}^{\infty}\frac{a_{2n}}{u^{2n}},\quad u{\to}\infty,\label{3.16}\end{equation}
where\\
$$a_{0}^{2}=1/4,\, 2a_{0}a_{2}=-1/4,\, a_{4}+a_{2}^{2}=3c_{2},\,
 a_{6}+ 2a_{2}a_{4}=5c_{4}$$\\
3)\begin{equation}|q(2\zeta)|^2{\sim}\pi(\frac{1}{2}u+\sum_{n=1}^{\infty}\frac{2nc_{2n}}{u^{2n+1}}),\quad u{\to}\infty.\label{3.17}\end{equation}
4)\begin{equation}q^{2}(2\zeta)=e^{iu}\pi[\sum_{n=0}^{\infty}\frac{a_{2n}(1-2n)}{u^{2n}}+
i\sum_{n=0}^{\infty}\frac{a_{2n}}{u^{2n-1}}],\quad u{\to}\infty.\label{3.18}\end{equation}\end{Tm}
\emph{Proof.}  We use the well-known system \cite{Meh},\cite{TrWi}:
\begin{equation}\frac{d}{d\zeta}[{\zeta}Q_{\zeta}(\zeta,\zeta)]=|r^{2}(\zeta)|,\quad 2\pi{\zeta}Q_{\zeta}(-\zeta,\zeta)=\Im[r^{2}(\zeta)]
\label{3.19}\end{equation}
\begin{equation}\frac{d}{d\zeta}[Q_{\zeta}(\zeta,\zeta)]=2Q_{\zeta}^{2}(-\zeta,\zeta)\quad
\frac{d}{d\zeta}[{\zeta}Q_{\zeta}(-\zeta,\zeta)]=\Re[r^{2}(\zeta)].
\label{3.20}\end{equation} We  need also the asymptotic relation (see \cite{TrWi}, formula (89)):
\begin{equation}-2{\zeta}Q_{\zeta}(\zeta,\zeta){\sim}-\frac{1}{4}u^{2}-
\frac{1}{4}+\sum_{n=1}^{\infty}\frac{c_{2n}}{u^{2n}},\quad u{\to}\infty,\quad (u=2\pi\zeta).\label{3.21}\end{equation}
where $c_2=-\frac{1}{4},\,c_4=-\frac{5}{2}.$ Formulas \eqref{3.19} and \eqref{3.21} imply
that
\begin{equation}|r(\zeta)|^{2}{\sim}\pi(\frac{1}{2}u+\sum_{n=1}^{\infty}\frac{2nc_{2n}}{u^{2n+1}}),\quad u{\to}\infty.\label{3.22}\end{equation}According to \eqref{3.21}we have
\begin{equation}Q_{\zeta}(\zeta,\zeta){\sim}\pi(\frac{1}{4}u+
\frac{1}{4u}-\sum_{n=1}^{\infty}\frac{c_{2n}}{u^{2n+1}}),\quad u{\to}\infty.\label{3.23}\end{equation}
Using relations \eqref{3.20} and (3.23) we derive that
\begin{equation}Q_{\zeta}^{2}(-\zeta,\zeta){\sim}\pi^{2}(\frac{1}{4}
-\frac{1}{4u^{2}}+\sum_{n=1}^{\infty}\frac{c_{2n}(2n+1)}{u^{2n+2}}),\quad u{\to}\infty.\label{3.24}\end{equation}
Consequently,
\begin{equation}Q_{\zeta}(-\zeta,\zeta){\sim}
\pi\sum_{n=0}^{\infty}\frac{a_{2n}}{u^{2n}},\quad u{\to}\infty,\label{3.25}\end{equation}
where in view of \eqref{3.24} we get
\begin{equation}a_{0}^{2}=1/4,\, 2a_{0}a_{2}=-1/4,\, 2a_{0}a_{4}+a_{2}^{2}=3c_{2},\,
 2a_{0}a_{6}+ 2a_{2}a_{4}=5c_{4}, ...\label{3.26}\end{equation}
Taking into account \eqref{3.19} and \eqref{3.20} we have
\begin{equation}r^{2}(\zeta)=
\frac{d}{d\zeta}[{\zeta}Q_{\zeta}(-\zeta,\zeta)]+i2\pi{\zeta}Q_{\zeta}(-\zeta,\zeta).\label{3.27}\end{equation}
According to \eqref{3.25} and \eqref{3.27} the equality
\begin{equation}r^{2}(\zeta){\sim}\pi\sum_{n=0}^{\infty}\frac{a_{2n}(1-2n)}{u^{2n}}+
i\pi\sum_{n=0}^{\infty}\frac{a_{2n}}{u^{2n-1}},\quad u{\to}\infty.\label{3.28}\end{equation}
holds.
Comparing formulas \eqref{3.5}, \eqref{3.10} with \eqref{3.22}, \eqref{3.23} and\eqref{3.25}, \eqref{3.28} we obtain the assertion of the theorem.
\begin{Rk}\label{Remark 3.4} In the next section we shall prove that
\begin{equation}a_{0}=1/2.\label{3.29}\end{equation}
Hence, using relations \eqref{3.26} we can find all coefficients $a_{2n}, \, (n=1,2,...).$\end{Rk}
In order to derive the asymptotic of $H_{2}(\zeta)$ we use the well-known Krein's formula (see \cite{GoKr}, Ch.IV):
\begin{equation}q_{1}^{2}(\zeta)=exp[2\int_{0}^{\zeta}R_{t}(t,0)dt].\label{3.30}\end{equation}
From \eqref{3.25}, \eqref{3.27} and \eqref{3.30} we deduce that
\begin{equation}\log[q_{1}^{2}(\zeta)]={\pi}\zeta+\beta+o(1),\quad \zeta{\to}\infty,\label{3.31}\end{equation}
where
\begin{equation}\beta=\int_{0}^{\zeta}[2R_{t}(t,0)-\pi]dt. \label{3.32}\end{equation}
It follows from \eqref{2.10} and \eqref{3.31} that
\begin{equation}\log[q_{2}^{2}(\zeta)]=-{\pi}\zeta-\beta-2\log{2}+o(1),\quad \zeta{\to}\infty.\label{3.33}\end{equation}
\begin{Rk}\label{Remark 3.5} Taking into account \eqref{2.6} and \eqref{3.17}, \eqref{3.18}, \eqref{3.29}
we obtain the asymptotic of $H_{1}(\zeta)$.\end{Rk}
\begin{Rk}\label{Remark 3.6} Taking into account \eqref{2.9} and \eqref{3.31}- \eqref{3.33}
we obtain the asymptotic of $H_{2}(\zeta)$.\end{Rk}
\emph{Let us investigate the expressions:}
$(S_{\zeta}^{-1}e^{it\pi},e^{it\pi})$ and $(S_{\zeta}^{-1}1,1),$ \emph{which are important in the theory of the integral operators with difference kernels}(see \cite{Sak4}).\\
According to  \eqref{2.7} and \eqref{2.10} we  have
\begin{equation}(S_{\zeta}^{-1}e^{it\pi},e^{it\pi})=\int_{0}^{\zeta}|q^{2}(t)|dt,\quad
(S_{\zeta}^{-1}1,1)=\int_{0}^{\zeta}|q_{1}^{2}(t)|dt\label{3.34}\end{equation}
Relations \eqref{3.17} and \eqref{3.34} imply that
\begin{equation}(S_{\zeta}^{-1}e^{it\pi},e^{it\pi})=\frac{\pi^{2}\zeta^{2}}{4}+O(1),\quad \zeta{\to}\infty.\label{3.35}\end{equation}
Relations \eqref{3.31} and \eqref{3.34} imply that
\begin{equation}(S_{\zeta}^{-1}1,1)=\frac{1}{\pi}e^{\pi\zeta+\beta}[1+o(1)],\quad \zeta{\to}\infty.\label{3.36}\end{equation}
\section{Fredholm determinants, integral representations}
Let us introduce the operators :
\begin{equation}K_{\pm}(\zeta)f=\int_{-1}^{1}k_{\pm}(x,t,\zeta)f(t)dt, \quad f(x){\in}L^{2}(-1,1),
\label{4.1}\end{equation}and
\begin{equation}K(\zeta)f=\int_{-1}^{1}k(x,t,\zeta)f(t)dt, \quad f(x){\in}L^{2}(-1,1),
\label{4.2}\end{equation}
where \begin{equation}k(x,t,\zeta)=\frac{\sin \zeta\pi(x-t)}{\pi(x-t)},\label{4.3}\end{equation} \begin{equation} k_{\pm}(x,t,\zeta)=[\frac{\sin \zeta\pi(x-t)}{\pi(x-t)}{\pm}\frac{\sin \zeta\pi(x+t)}{\pi(x+t)}]/2
\label{4.4}\end{equation}In this section we shall investigate the operators
\begin{equation}S(\zeta,\lambda)=I-{\lambda}K(\zeta)\quad and \quad S_{\pm}(\zeta,\lambda)=I-{\lambda}K_{\pm}(\zeta),\quad 0<\lambda{\leq}1.
\label{4.5}\end{equation}
It is easy to see that
\begin{equation} K_{\pm}=\frac{I{\pm}J}{2}K=K\frac{I{\pm}J}{2},\label{4.6}\end{equation}
where $Jf(x)=f(-x)$.

\begin{La}\label{Lemma 4.1}The following equality holds (see\cite{TrWi}):
\begin{equation}(I-{\lambda}K_{\pm})^{-1}K_{\pm}=\frac{I{\pm}J}{2}(I-{\lambda}K)^{-1}K\label{4.7}\end{equation}\end{La}
\emph{Proof.} It is easy to see that
\begin{equation}(\frac{I{\pm}J}{2})^{n}=\frac{I{\pm}J}{2},\quad n=1,2,..\label{4.8}\end{equation}
Using relations \eqref{4.6},  \eqref{4.8} and inequality $\norm{K}<1$ we obtain
\begin{equation}(I-{\lambda}K_{\pm})^{-1}-I=\frac{I{\pm}J}{2}\sum_{n=1}^{\infty}K^{n}\lambda^{n}=
\frac{I{\pm}J}{2}[(I-{\lambda}K)^{-1}-I]\label{4.9}\end{equation}The assertion of Lemma 4.1 follows from equality
 \eqref{4.9}.\\
We consider the Fredholm determinants
 \begin{equation}P_{\pm}(\zeta,\lambda)=\det(I-{\lambda}K_{\pm}(\zeta)),\quad P(\zeta,\lambda)=\det(I-{\lambda}K(\zeta))
\label{4.10}\end{equation}
\begin{La}\label{Lemma 4.2} The following relation is valid
\begin{equation}\frac{d}{d\zeta}\log{P_{\pm}(\zeta,\lambda)}=-\lambda\{((I-{\lambda}K(\zeta))^{-1}\phi_{1},\phi_{1})
{\pm}\Re[((I-{\lambda}K(\zeta))^{-1}\phi_{1},\overline{\phi_{1}})]\},\nonumber\end{equation}
where $\phi_{1}(x,\zeta)=e^{ix\pi\zeta}/\sqrt{2}.$\end{La}
\emph{Proof.} Let us write the equality
\begin{equation}\frac{d}{d\zeta}\log{P_{\pm}(\zeta,{\lambda})}=-{\lambda}tr[(I-{\lambda}K_{\pm}(\zeta))^{-1}\frac{d}{d\zeta}K_{\pm}(\zeta)].
\label{4.11}\end{equation}
It follows from \eqref{4.4}, \eqref{4.8} and \eqref{4.11} that
\begin{equation}\frac{d}{d\zeta}\log{P_{\pm}(\zeta,\lambda)}=-{\lambda}tr[(I-{\lambda}K(\zeta))^{-1}\frac{I{\pm}J}{2}\frac{d}{d\zeta}K(\zeta)].
\label{4.12}\end{equation}
According to \eqref{4.2} we have
\begin{equation}\frac{d}{d\zeta}K(\zeta)f=\frac{1}{2}\int_{-1}^{1}
[e^{i\pi\zeta(x-t)}+e^{-i\pi\zeta(x-t)}]f(t)dt
\label{4.13}\end{equation}We introduce the one-dimensional operators
\begin{equation}T_{\pm}(\zeta,\lambda)f=\frac{\lambda}{4}(I-{\lambda}K(\zeta))^{-1}\int_{-1}^{1}
e^{{\pm}i\pi\zeta(x-t)}f(t)dt, \label{4.14}\end{equation}
\begin{equation}V_{\pm}(\zeta,\lambda)f=\frac{\lambda}{4}(I-{\lambda}K(\zeta))^{-1}\int_{-1}^{1}
e^{{\pm}i\pi\zeta(-x-t)}f(t)dt, \label{4.15}\end{equation}
Let us consider a complete orthonormal system functions $\phi_{n}(x,\zeta), \quad (n=1,2,...)$ in the space
$L^2(-1,1)$ such, that $\phi_{1}(x,\zeta)=e^{ix\pi\zeta}/\sqrt{2}.$ Then we have
\begin{equation}tr{T_{+}(\zeta,\lambda)}=\frac{\lambda}{4}\sum_{n=1}^{\infty}((I-{\lambda}K(\zeta))^{-1}\phi,\phi_{n})(\phi,\phi_{n}),\quad
\phi=\sqrt{2}\phi_{1}.
\label{4.16}\end{equation}It follows from \eqref{4.16}, that
\begin{equation}tr{T_{+}(\zeta,\lambda)}=\frac{\lambda}{2}((I-{\lambda}K(\zeta))^{-1}\phi_{1},\phi_{1})
\label{4.17}\end{equation}
In the same way we obtain the relations
\begin{equation}tr{T_{-}(\zeta,\lambda)}=
\frac{\lambda}{2}((I-{\lambda}K(\zeta))^{-1}\overline{\phi_{1}},\overline{\phi_{1})},
\label{4.18}\end{equation}
\begin{equation}tr{V_{+}(\zeta,\lambda)}=
\frac{\lambda}{2}((I-{\lambda}K(\zeta))^{-1}\overline{\phi_{1}},\phi_{1}),\label{4.19}\end{equation}
\begin{equation}tr{V_{-}(\zeta,\lambda)}=
\frac{\lambda}{2}((I-{\lambda}K(\zeta))^{-1}\phi_{1},\overline{\phi_{1})},
\label{4.20}\end{equation}The kernel of the operator $K$ is real. Hence, the relations \eqref{4.17}-\eqref{4.20} imply, that
\begin{equation}tr{T_{-}(\zeta,\lambda)}=tr{T_{+}(\zeta,\lambda)},\quad tr{V_{-}(\zeta,\lambda)}=\overline{trV_{+}(\zeta,\lambda)}.\label{4.21}\end{equation}
Taking into account \eqref{4.22} we deduce
\begin{equation}\frac{d}{d\zeta}\log{P_{\pm}(\zeta,\lambda)}=-\lambda\{tr{T_{+}(\zeta,\lambda)}{\pm}\Re[tr{V_{-}(\zeta,\lambda)]}\}.\label{4.22}\end{equation}
The assertion of the lemma follows directly from relations
\eqref{4.17}, \eqref{4.19} and \eqref{4.22}.\\
Further we need the operator $Uf(x)=g(s)=f(s/\zeta)/\sqrt{\zeta}$, which maps isometrically $L^{2}(-1,1)$
onto $L^{2}(-\zeta,\zeta)$. It is easy to see that  $U^{-1}g(s)=f(x)=g(x\zeta)\sqrt{\zeta}$ and $U^{*}=U^{-1}.$ It follows from \eqref{4.2} and \eqref{4.3} that
\begin{equation}K(\zeta)=U^{-1}C_{\zeta}U,\quad K_{\pm}(\zeta)=U^{-1}C_{\pm}(\zeta)U\label{4.23}\end{equation}
where
\begin{equation}C_{\zeta}g(s)=
\int_{-\zeta}^{\zeta}
\frac{\sin 
\pi(s-t)}{\pi(s-t)}g(s)ds\label{4.24}\end{equation}
\begin{equation}C_{\pm}(\zeta)g(s)=\frac{1}{2}
\int_{-\zeta}^{\zeta}
[\frac{\sin 
\pi(s-t)}{\pi(s-t)}{\pm}\frac{\sin 
\pi(s+t)}{\pi(s+t)}]g(s)ds\label{4.25}\end{equation}
Using relations \eqref{4.10} and \eqref{4.23} we obtain, that
 \begin{equation}P_{\pm}(\zeta,\lambda)=\det(I-{\lambda}C_{\pm}(\zeta))^{-1},\quad P(\zeta,\lambda)=\det(I-{\lambda}C_{\zeta})^{-1}
\label{4.26}\end{equation}
Lemma 4.2 and relations \eqref{4.23}, \eqref{4.26} imply the assertion:
\begin{La}\label{Lemma 4.3} The following relation is valid
\begin{equation}2\zeta\frac{d}{d\zeta}\log{P_{\pm}(\zeta,\lambda)}=-\lambda\{((I-{\lambda}C_{\zeta})^{-1}\psi_,\psi)
{\pm}\Re[((I-{\lambda}C_{\zeta})^{-1}\psi,\overline{\psi})]\},
\nonumber\end{equation} where $\psi(x)=e^{ix\pi}$.
\end{La}

 Using relations \eqref{3.2} and Lemma 4.3 we get
\begin{La}\label{Lemma 4.4}If $\lambda=1$, then
\begin{equation}2\zeta\frac{d}{d\zeta}\log{P_{\pm}(\zeta,1)}=-\{(S_{2\zeta}^{-1}\psi_{1},\psi_{1
})
{\pm}\Re[(S_{2\zeta}^{-1}\psi_{1},\overline{\psi_{1}})]\},\label{4.27}\end{equation}
where $\psi_{1}(x)=U_{\zeta}\psi(x)=\psi(x-\zeta).$\end{La}
It follows from \eqref{2.7} and \eqref{4.27} that
\begin{equation}2\zeta\frac{d}{d\zeta}\log{P_{\pm}(\zeta,1)}=-\{\int_{0}^{2\zeta}|q^{2}(s)|ds{\pm}
\Re[\int_{0}^{2\zeta}e^{2is\pi}q^{2}(s)ds]\}.\label{4.28}\end{equation}
Taking into account \eqref{3.10}, \eqref{3.11} and \eqref{4.28}
we obtain
\begin{Tm}\label{Theorem 4.5}If $\lambda=1$, then
\begin{equation}\zeta\frac{d}{d\zeta}\log{P_{\pm}(\zeta,1)}=-\{\int_{0}^{\zeta}|r^{2}(s)|ds{\pm}
\Re[\int_{0}^{\zeta}r^{2}(s)ds]\}.\label{4.29}\end{equation}\end{Tm}
Now we show that Remark 3.4 is valid.
\begin{Cy}\label{Corollary 4.6}The equality
\begin{equation}a_{0}=1/2 \label{4.30}\end{equation}
holds.\end{Cy}
\emph{Proof.} According to \eqref{3.28} we have
\begin{equation}\Re[r^2(\zeta)]={\pi}a_{0}+o(1),\quad \zeta{\to}\infty.\label{4.31}\end{equation}
Taking into account \eqref{4.29} and \eqref{4.31} we obtain
\begin{equation}\frac{d}{d\zeta}[\log{P_{+}(\zeta,1)}-\log{P_{-}(\zeta,1)}]=-
2\pi[{a_{0}}+o(1)].\label{4.32}
\end{equation}
Hence, the following relation is valid
\begin{equation}\log{P_{+}(\zeta,1)}-\log{P_{-}(\zeta,1)}=-
2\pi\zeta[{a_{0}}+o(1)].\label{4.33}
\end{equation} Using  asymptotic formula (90) from the paper \cite{TrWi} we have
\begin{equation}\log{P_{+}(\zeta,1)}-\log{P_{-}(\zeta,1)}=-
\pi\zeta[1+o(1)].\label{4.34}
\end{equation}Relations \eqref{4.33} and \eqref{4.34} imply relation \eqref{4.30}.
The Corollary is proved.\\
Let us introduce the notations
\begin{equation}\sigma_{\pm}(\zeta)=\zeta\frac{d}{d\zeta}\log{P_{\pm}(\zeta,1)},\quad
\sigma(\zeta)=\zeta\frac{d}{d\zeta}\log{P(\zeta,1)}.
\label{4.35}\end{equation} In view of \eqref{4.7} we get
\begin{equation}\sigma_{+}(\zeta)+\sigma_{-}(\zeta)=\sigma(\zeta).\label{4.36}\end{equation}
Theorem 4.5 and relation \eqref{4.36} imply the assertion:
\begin{Cy}\label{Corollarry 4.7}If $\lambda=1$ then
\begin{equation}\zeta\frac{d}{d\zeta}\log{P(\zeta,1)}=-\int_{0}^{\zeta}|r^{2}(s)|ds.\label{4.37}
\end{equation}\end{Cy}
We note that relations \eqref{4.36} and \eqref{4.37} are well-known \cite{TrWi}.
\begin{Rk}\label{Remark 4.8}Theorem 4.5 (see formula\eqref{4.29}) can be used by solving the following problems:\\
1. Find the asymptotics of $\zeta\frac{d}{d\zeta}\log{P_{\pm}(\zeta,1)}$ when $\zeta{\to}\infty$ with the help of formulas $\eqref{3.22}$ and $\eqref{3.28}$.\\
2. Find the asymptotics of $\zeta\frac{d}{d\zeta}\log{P_{\pm}(\zeta,1)}$ when $\zeta{\to}0$ with the help of formula $\eqref{3.11}$.\\
3. Find the estimation of $\zeta\frac{d}{d\zeta}\log{P_{\pm}(\zeta,1)}$ with the help of formula $\eqref{3.11}.$\end{Rk}

\end{document}